\documentclass[a4paper,12pt]{amsart}

\usepackage{latexsym}
\usepackage{amssymb}
\usepackage{amsthm}
\usepackage{amscd}
\usepackage{amsmath}
\usepackage{graphics}
\usepackage[all]{xy}

\newtheorem{Theorem}{Theorem}[section]
\newtheorem{Lemma}[Theorem]{Lemma}
\newtheorem{Corollary}[Theorem]{Corollary}
\newtheorem{Proposition}[Theorem]{Proposition}
\newtheorem{Remark}[Theorem]{Remark}

\newtheorem{Example}[Theorem]{Example}

\newtheorem{Definition}[Theorem]{Definition}

\def\fkm{{\frak m}}

\def\fkn{{\frak n}}

\def\opn#1#2{\def#1{\operatorname{#2}}}

\opn\Spec{Spec}
\opn\Supp{Supp}
\opn\supp{supp}
\opn\Max{Max}
\opn\max{max}
\opn\Min{Min}
\opn\min{min}
\opn\Ass{Ass}
\opn\Assh{Assh}
\opn\depth{depth}
\opn\rank{rank}

\opn\Sym{Sym}
\def\Rees{{\mathcal R}}

\opn\div{div}
\opn\Div{Div}
\opn\cl{cl}
\opn\Cl{Cl}

\opn\Ker{Ker}
\opn\Coker{Coker}
\opn\Im{Im}
\opn\Hom{Hom}
\opn\Tor{Tor}
\opn\Ext{Ext}
\opn\End{End}
\opn\Aut{Aut}
\opn\id{id}

\opn\nat{nat}
\opn\pff{pf}
\opn\Pf{Pf}
\opn\GL{GL}
\opn\SL{SL}
\opn\G{G}
\opn\E{E}
\opn\H{H}
\opn\M{M}
\opn\mod{mod}
\opn\ord{ord}
\opn\det{det}
\opn\I{I}    
\opn\Soc{Soc}
\opn\Fitt{Fitt}

\opn\chara{char}
\opn\length{\ell}
\opn\pd{pd}
\opn\rk{rk}
\opn\projdim{proj\,dim}
\opn\injdim{inj\,dim}
\opn\rank{rank}
\opn\depth{depth}
\opn\grade{grade}
\opn\height{ht}
\opn\embdim{emb\,dim}
\opn\codim{codim}

\renewcommand{\tilde}{\widetilde}
\renewcommand{\bar}{\overline}
\renewcommand{\hat}{\widehat}

\title{Modules of reduction number one}

\author{Futoshi Hayasaka}

\dedicatory{Dedicated to Professor Shiro Goto on the occasion of his sixtieth birthday}

\address{Department of Mathematics, School of Science and Technology,
Meiji University, 1-1-1, Higashimita, Tama-ku, Kawasaki, Kanagawa, 
214--8571, JAPAN}

\email{fhayasaka@math.meiji.ac.jp}

\date{}


\begin{document}
\maketitle

\setlength{\baselineskip}{16pt}

\begin{abstract}
Let $(A, \fkm)$ be a Noetherian local ring and $N$ a parameter module in $F=A^r$ and $M=N:_F \fkm$ the socle module of $N$. In this paper, we shall prove that the module $M=N:_F \fkm$ has a reduction number at most one and hence its Rees algebra 
$\Rees(M)$ is Cohen-Macaulay, if the base ring $A$ is Cohen-Macaulay of dimension two and the rank of $N$ is greater than or equal to two. This result gives numerous examples of Cohen-Macaulay Rees algebras of modules, which are not integrally closed and not a parameter module. 
\end{abstract}


\section*{Introduction}

Let $(A, \fkm)$ be a Noetherian local ring of dimension $d=\dim A$ and 
$F=A^r$ a free module of rank $r>0$. 
Let $M$ be a submodule of $F$ such that the length $\ell_A(F/M)$ of the quotient $F/M$ is finite. 
Then the Rees algebra of the module $M$ is defined to be 
the image of the natural homomorphism from $\Sym_A(M)$ to $\Sym_A(F)$, 
which is a subalgebra of the polynomial ring over $A$ (cf. \cite{EHU, SUV}). 
We denote by $\Rees(M)$ the Rees algebra of $M$. 

In 1997, Katz and Kodiyalam showed in \cite{KK} that the Rees algebra 
$\Rees(M)$ is Cohen-Macaulay if and 
only if the reduction number of $M$ is at most one when $A$ is 
Cohen-Macaulay of dimension two (see also \cite{SUV}). When $r=1$, this can be 
viewed as the module version of the Goto-Shimoda theorem on the 
Cohen-Macaulayness of the Rees algebra of ideals in dimension two (\cite{GS}). 
In the same paper \cite{KK}, they 
also showed that 
the reduction number of complete (i.e. integrally closed) modules 
over a two dimensional regular local ring is at most one, and hence 
its Rees algebra is Cohen-Macaulay. 
When $r=1$, this can be also viewed as a generalization of the result of Huneke and Sally \cite{HS} on integrally closed ideals in a regular local ring of dimension two. 

Inspired by their results, to construct a large class of modules 
whose reduction number is at most one, we shall study the socle module of a parameter module of any rank. When $r=1$, this problem was first 
investigated by Corso, Polini and Vasconcelos in \cite{CPV, CP}. 
When $A$ is Cohen-Macaulay, Corso, Huneke, Polini and Vasconcelos 
gave a complete answer to this problem in their papers \cite{CHV,CPV,CP}. 
Let $Q$ be a parameter ideal in $A$ and $I=Q:\fkm$ the socle ideal of $Q$. 
Then Corso and Polini showed in \cite{CP} that
 the equality $I^2=QI$ holds if $A$ is Cohen-Macaulay but not regular. 
In \cite{CHV}, Corso, Huneke and Vasconcelos showed that
the socle ideal of an $\fkm$-primary Gorenstein ideal contained in $\fkm^2$
has a reduction number at most one if the embedding dimension of $A$ 
is at least two. The problem is related to the integral closedness of parameter ideals. In \cite{G}, Goto studied certain integrally closed ideals and showed that $A$ is regular if there exists an integrally closed parameter ideal $Q$ in $A$. He also showed that $Q$ has a system of generators of special type in this case. 
So, we know many things in the case where $r=1$ and $A$ is Cohen-Macaulay. 

The purpose of this paper is to give a natural generalization of their results for modules and construct a large class of modules whose reduction 
number is at most one. In particular, our focus 
is in the case where the rank $r$ is at least two. The main result is 
the following.

\setcounter{section}{5}

\begin{Theorem}
Let $(A, \fkm)$ be a Cohen-Macaulay local ring of dimension two 
and $N$ a parameter module in a free module $F=A^r$ of rank $r>0$, 
namely, $\ell_A(F/N) < \infty$ and $\mu_A(N)=r+1$. 
Assume that $N \subseteq \fkm F$ and put $M=N:_F \fkm$ the socle module of $N$. Then the following three conditions are equivalent. 
\begin{enumerate}
\item[$(1)$] $M^2 \ne NM$. 
\item[$(2)$] $N$ is an integrally closed module in $F$. 
\item[$(3)$] $A$ is a regular local ring, $r=1$ and $N$ is an integrally closed ideal in $A$. 
\end{enumerate}
\end{Theorem}

\noindent
Here we consider the products and the powers of modules 
inside the symmetric algebra $\Sym_A(F)$. 

\hspace{0.5mm}

This result is closely related to the result of Simis-Ulrich-Vasconcelos. 
They showed in \cite{SUV} that the
Rees algebras of modules of this type are always Cohen-Macaulay in more 
general situation. But our method of proof is different from theirs, 
and the heart of our theorem is in the condition of rank one. 
When $r=1$, the above theorem is already known. Consequently, what we have to prove is that the rank $r$ is just one if $M^2 \ne NM$ holds or $N$ is an integrally closed in $F$. 

We explain the construction of this paper. 
In section 1, we fix our notations and definitions. 
In section 2, we introduce a notion of a perfect matrix, which is a general
notion of a regular sequence in some sense. 
We shall prove that modules generated by a perfect matrix have the same 
property on ideals generated by a regular sequence. 
In section 3, we shall give a generalization of results in \cite{GH2}
and discuss integrally closed modules with finite homological dimension. 
We shall use these results in the following sections and the proof 
of the main theorem. Also, we shall discuss an analogue of the Levin-Vasconcelos theorem in \cite{LV} and give a generalization of the result in \cite{AP}. 
In section 4, we shall give two classes of modules with reduction number one, 
which is a natural generalization of the results in \cite{CHV, CPV, CP}. 
Finally, in section 5, we shall give a proof of Theorem 5.1 and give 
some examples of modules whose reduction number is at most one. 

\setcounter{section}{0}

\section{Notation}

In this section, we shall establish our notations and definitions. 
Let $A$ be a commutative Noetherian ring and $F$ a free $A$-module 
of rank $r>0$. Let $M$ be a submodule of $F$. 
Then the Rees algebra of the module $M$ is defined to be 
the image of the natural homomorphism from $\Sym_A(M)$ to $\Sym_A(F)$, 
which is a subalgebra of the polynomial ring over $A$. 
We denote by $\Rees(M)$ the Rees algebra of $M$. 
We always consider a module $M$ in a fixed free module $F$. 
We refer the reader to \cite{EHU} for the other 
definition of the Rees algebras of 
modules, which does not depend on the embedding of the module in a free module. 
Let $S=\Sym_A(F)$ be the symmetric algebra of a fixed free module $F$, which is a polynomial ring over $A$. Then the module $F$ can be identified with the set of linear forms $S_1$ in $S$. An element $f \in F$ is said to be integral over $M$ if the linear form $f \in S_1$ is integral over the Rees algebra $\Rees(M)$ of $M$. The set of elements in $F$ which are integral over $M$ is called the integral closure of $M$ in $F$, and is denoted by $\bar{M}$. The integral closure of $M$ is a submodule of $F$ containing $M$. $M$ is said to be integrally closed in $F$ if the equality $\bar{M}=M$ holds. 
For each integer $n \geq 0$, let $M^n$ denote the homogeneous component of 
degree $n$ in $\Rees(M)$. 
For another submodule $N$ of $F$, let 
$NM$ denote the product of two modules $N$ and $M$. 
Here we always consider the products and the powers of modules 
inside the symmetric algebra $S=\Sym_A(F)$ of a fixed free module $F$. 
We say that $N$ is a reduction of $M$ if $N$ is a submodule of $M$ and the 
equality $M^{n+1}=NM^n$ holds in $S_{n+1}$ for some integer $n \geq 0$. 
Equivalently, $N$ is a reduction of $M$ if and only if the ring extension 
$\Rees(N) \subseteq \Rees(M)$ is integral. Also, a submodule $N$ of $M$ is a reduction of $M$ if and only if $M \subseteq \bar{N}$. 
A reduction of $M$ is said to be minimal, if it has no proper reductions. 
Furthermore, we assume that $A$ is a local ring with the maximal ideal $\fkm$ of dimension $d$. Then, for any reduction $N$ of $M$, one has the inequality $\mu_A(N) \geq \lambda(M)$ where $\mu_A(-)$ denotes the minimal number of generators of the inside module and $\lambda(M)=\dim \Rees(M)/\fkm \Rees(M)$ the analytic spread of $M$. A reduction $N$ of $M$ is minimal if the equality $\mu_A(N) = \lambda(M)$ holds. Note that the analytic spread $\lambda(M) = d+r-1$, if the colength $\ell_A(F/M)$ is finite. Following \cite{KK}, a submodule $N$ of $F$ is said to be a parameter module in $F$, if the colength $\ell_A(F/N)$ is finite and the equality $\mu_A(N)=d+r-1$ holds. For further facts and details on the Rees algebras of modules, we refer the reader to \cite{SUV}. 

\section{Perfect matrix}

In this section, we shall introduce a notion of a perfect matrix and 
investigate the property of modules with a perfect matrix. 

Let $A$ be a commutative Noetherian ring and $n \geq r > 0$ 
integers. 
Let $F=A^r$ be a free module over $A$ with rank $r>0$. 
We fix a free basis $ \{ t_1, \dots , t_r \} $ for $F$. Let 
$S=\Sym_A(F)$ be the symmetric algebra of $F$, which is 
identified with the polynomial ring 
$A[t_1, \dots , t_r]$ over $A$. 
Fixing a free basis for $F$, every element $f$ of $F$ can be written 
as a column with $r$ entries $f=^t(f_1, \dots , f_r)$. 
Furthermore, since $F$ can be canonically identified with $S_1$, 
we can write $f \in F$ as a linear form $f_1 t_1+ \dots + f_r t_r$ in $S$. 
Let $N$ be a submodule of $F$ generated by 
$c_1, \dots , c_n$. Fixing a free basis for $F$, 
we have the matrix $\tilde{N}=(c_{ij})$ associated to the module $N$, where 
$ c_j=c_{1j}t_1 + \dots + c_{rj}t_r,  \  (c_{ij} \in A). $
Then we first introduce a notion of a perfect matrix. 

\begin{Definition}
A matrix $\varphi$ of size $r \times n$ over $A$ is called 
perfect if the following two conditions are satisfied: 
$(i)$ $I_r(\varphi) \subsetneq A$, and $(ii)$ $I_r(\varphi)$ has the greatest 
possible grade $n-r+1$. 
Here we denote by $I_r(\varphi)$ the ideal generated by $r \times r$-minors 
of $\varphi$. 
\end{Definition}

This is a general notion of a regular sequence in some sense. 

\begin{Example} The following matrices are perfect. 
\begin{enumerate}
\item[$(1)$] A generic matrix $(X_{ij})$. 
\item[$(2)$] A matrix $\tilde{N}$ of a parameter module $N$ over 
a Cohen-Macaulay local ring. {\rm{(}}The matrix $\tilde{N}$ is a parameter 
matrix in the sense of Buchsbaum-Rim \cite{BR}. {\rm{)}} 
\item[$(3)$] A matrix 
$\left(
\begin{array}{cccccc}
a_1 & \cdots & a_{\ell} & 0        & \cdots &  0        \\
0   & a_1    & \cdots   & a_{\ell} & \ddots &  \vdots   \\
\vdots & \ddots& \ddots &          & \ddots &  0        \\
0   & \cdots & 0        & a_1      & \cdots &  a_{\ell} 
\end{array}
\right), $
where $a_1, \dots , a_{\ell}$ is a regular sequence on $A$. 
\end{enumerate}
\end{Example}

A module with a perfect matrix has the following property, which is a generalization of the well-known fact in the case where the ideal generated by a regular sequence. 

\begin{Proposition}
Let $N$ be a submodule of $F$ generated by 
$c_1, \dots , c_n$. Assume that the matrix $\tilde{N}$ of $N$ is perfect and 
$n>r$. If $h_1, \dots , h_n \in F$ and satisfy the relation 
$$h_1c_1 + \dots + h_nc_n=0 \ \ \mbox{in} \ \ S_2, $$
then we have 
$$h_i \in (c_1, \dots , \hat{c_i}, \dots , c_n) \ \ \ \mbox{for all} \ \ 
1 \leq i \leq n. $$
\end{Proposition}

\begin{proof}
Since $n>r$ and $\grade I_r(\tilde{N})=n-r+1$, 
the generalized Koszul complex (due to D. Kirby \cite{K}) 
$K_{\bullet}(\tilde{N}; 2)$ 
associated to a matrix $\tilde{N}$ and an integer $2$ is acyclic (cf. \cite{K} and \cite[Appendix A2.6]{E}). In particular, the first homology module of this complex 
$H_1(K_{\bullet} (\tilde{N}; 2))$ vanishes.
By a construction of the generalized Koszul complex, $H_1(K_{\bullet} (\tilde{N}; 2)) \cong H_1(K_{\bullet} (c_1, \dots , c_n ; S))_2$, where $K_{\bullet} (c_1, \dots , c_n ; S)$ is the ordinary graded Koszul complex associated to a sequence $c_1, \dots , c_n$ in $S$. 
Looking at the graded component of degree $2$ of the Koszul complex 
$K_{\bullet} (c_1, \dots , c_n ; S)$, we have the exact sequence 
$$
\xymatrix{
\wedge^2 \otimes S_0 \ar[r]^{d_2} & \wedge^1 \otimes S_1 \ar[r]^{d_1} & 
\wedge^0 \otimes S_2, 
}
$$
where $\wedge$ is the exterior algebra of a free module $G=A^n$ of rank $n$. 
Let $\{e_1, \dots ,e_n\}$ be a free basis for $G$. Then one can easily check that  

\noindent
{\bf{Claim 1}} \ 
$\displaystyle d_1 \left( \sum_{i=1}^n e_i \otimes h_i \right) 
=\sum_{i=1}^nh_i c_i $. 

\noindent
Since $\sum_{i=1}^nh_i c_i=0$, there exists an element 
$\xi \in \wedge^2 \otimes S_0$ such that 
$$d_2(\xi)=\sum_{i=1}^n e_i \otimes h_i. $$
Write 
$\displaystyle{\xi = \sum_{1 \leq i<j \leq n}(e_i \wedge e_j) \otimes a_{ij}}$, where 
$a_{ij} \in A$. Then one can also check that 

\noindent
{\bf{Claim 2}} \ 
$\displaystyle d_2(\xi)=\sum_{k=1}^n e_k \otimes \left( 
\sum_{1 \leq i<k \leq n} 
a_{ik}c_i - \sum_{1 \leq k<j \leq n} a_{kj}c_j \right). $

\noindent
Consequently, we have 

\begin{eqnarray*}
d_2(\xi) &=& \sum_{i=1}^n e_i \otimes h_i \\
    &=& \sum_{k=1}^ne_k \otimes \left( \sum_{1 \leq i<k \leq n} a_{ik}c_i - 
\sum_{1 \leq k<j \leq n} a_{kj}c_j \right). 
\end{eqnarray*}
By comparing the coefficient of each $e_i$, we get that  
$$h_i \in (c_1, \dots , \hat{c_i}, \dots , c_n) \ \ \ 
\mbox{for all} \ \ \ 1 \leq i \leq n. $$
\end{proof}

Proposition 2.3 is not true in the case where $n=r$ (i.e., 
$\tilde{N}$ is a square perfect matrix). There is the following 
simple example. 

\begin{Example}
Let $A=K[X]$ be a polynomial ring over a field
$K$. 
Let $N$ be a submodule of $F=A^2=Av+Aw$ such that the matrix 
$\tilde{N}$ is given by 
$$\left(
\begin{array}{ccccccc}
X & 0 \\ 
0 & X
\end{array}
\right),  $$
a square perfect matrix over $A$. Then we have the following relation:
$$w(Xv)-v(Xw)=0. $$ 
This means that Proposition 2.3 is not true in this case. 
\end{Example}

As a consequence of Proposition 2.3, we have the following. This can be viewed
as a part of the module version of Rees' theorem for ideals generated by 
a regular sequence \cite{Re}. We shall discuss in the forthcoming 
paper \cite{HH} for more results on this topic. 

\begin{Corollary}
Let $N$ be a submodule of $F$ generated by 
$c_1, \dots , c_n$. Assume that the matrix $\tilde{N}$ of $N$ is perfect 
and $n>r$. Then the following natural surjective homomorphism 
$$
\xymatrix{
(F/N)^{\oplus n} \ar@{>>}[r]^{[c_1 \cdots c_n]} & NF/N^2 
}
$$
is an isomorphism. 
\end{Corollary}

\begin{proof}
Let $h_i \in F$ and $\sum_{i=1}^nh_ic_i \in N^2$. Then 
there exist elements $g_i \in N$ such that 
$\sum_{i=1}^nh_ic_i=\sum_{i=1}^ng_ic_i$. Hence $\sum_{i=1}^n(h_i-g_i)c_i=0$. 
Since $\tilde{N}$ is perfect, each coefficient $h_i-g_i$ of $c_i$ is in $N$ by Proposition 2.3. Therefore each $h_i \in N$. Hence we have the isomorphism. 
\end{proof}

\section{Integrally closed modules with finite homological dimensions}

The main result in this section is the following. 

\begin{Theorem}
Let $(A,\fkm)$ be a Noetherian local ring and $M$ a finitely generated 
$A$-module such that $\depth_A M>0$. Assume that there exist two submodules 
$N$ and $L$ of $M$ such that 
\begin{enumerate}
	\item[$(1)$] $N \subseteq L \subseteq N:_M\fkm$,
	\item[$(2)$] $\fkm N \ne \fkm L$, and 
	\item[$(3)$] $\pd_A N<\infty$. 
\end{enumerate}
Then $A$ is a regular local ring. 
\end{Theorem}

\begin{proof}
By assumption $(2)$, we can take $xy \in \fkm L \backslash \fkm N$ 
such that $x \in \fkm \backslash \fkm^2$ and $y \in L$. 
Since $\depth_A N > 0$, we can choose $x$ as a non-zero-divisor on $N$. 
Note that $xy \in N$. 
Since $xy \notin \fkm N$, we can write $ N=Axy+\sum_{i = 1}^\ell Az_i$, where 
$\ell=\mu_A(N)-1$. 
Then we have the following. 

\smallskip

\noindent
{\bf{Claim}}
$N/xN=A \bar{xy} \oplus \sum_{i = 1}^\ell A \bar{z_i}$. 

\smallskip

\noindent
Indeed, let $a , b_i \in A$ and assume that 
$a(xy) + \sum_{i=1}^\ell {b_iz_i} \in xN$. 
Since $xN \subseteq \fkm N$, we have $a , b_i \in \fkm$. 
Hence $a(xy)=x(ay) \in xN$. We get the Claim. 

\noindent
Since $\fkm (xy) \subseteq xN$, $A \bar{xy} \cong A/\fkm$. 
Therefore we have $\pd_A A/\fkm < \infty$ because the residue field $A/\fkm$ is a direct summand of $N/xN$ with finite projective dimension. Consequently, 
$A$ is a regular local ring. 
\end{proof}

Theorem 3.1 is a slight general form of \cite[Theorem 1.1]{GH2} even if $N$ and $L$ are ideals in $M=A$. 
As a consequence, we have the following, which was first discovered by Burch \cite{Bu} when $r=1$ (cf. \cite{GH1, GH2}). 

\begin{Corollary}
Let $(A, \fkm)$ be a Noetherian local ring of positive depth and $F$ a free $A$-module of rank $r>0$. Assume that there exists a submodule $N$ of $F$ such that 
\begin{enumerate}
\item [$(1)$] $N$ is integrally closed in $F$, 
\item [$(2)$] $\fkm \in \Ass_A F/N$, and 
\item [$(3)$] $\pd_A N < \infty$. 
\end{enumerate}
Then $A$ is a regular local ring. 
\end{Corollary}

\begin{proof}
Assume $A$ is not regular. By Theorem 3.1, we have $\fkm (N:_F \fkm)=\fkm N$. Hence the module $N:_F \fkm $ is integral over $N$, using a standard technique on the determinant. By assumption $(1)$, we have the equality 
$N:_F\fkm=N$. This is contradict to assumption $(2)$. 
Hence $A$ is a regular local ring. 
\end{proof}

As the other consequence of Theorem 3.1, we have the following, which is the 
Levin-Vasconcelos theorem \cite{LV}. 

\begin{Corollary}
Let $(A, \fkm)$ be a Noetherian local ring and $M$ a finitely generated $A$-module with $\depth_A M>0$. If $\pd_A \fkm^nM < \infty$ for some positive integer $n>0$, then $A$ is a regular local ring. 
\end{Corollary}

\begin{proof}
Let $N=\fkm^nM$ and $L=\fkm^{n-1}M$. Then, by Nakayama's Lemma, $\fkm N \ne \fkm L$, because $N$ is non-zero. Hence $A$ is regular by Theorem 3.1. 
\end{proof}

\begin{Remark}{\rm{
All the results in this section work for the other homological dimensions 
(for instance, Gorenstein dimension \cite{Au}) 
instead of the projective dimension. 
In \cite{AP}, 
Asadollahi and Puthenpurakal showed Corollary 3.3 for various homological dimensions under the assumption $n \gg 0$. So our result Corollary 3.3 is a generalization of \cite{AP}. 
}}
\end{Remark}

\section{Classes of modules of reduction number one}

In this section, we shall give two classes of modules whose reduction number is at 
most one. 
Let $(A, \fkm)$ be a Noetherian local ring and $F$ a free module of rank $r>0$.
We begin with the following. 

\begin{Lemma}
Let $N$ be a submodule of $F$ generated by $c_1, \dots , c_n$, 
where $n=\mu_A(N)$. Assume $n>r$. Put $M=N:_F \fkm$ 
the socle module of $N$. If 
\begin{enumerate}
\item [$(1)$] $\tilde{N}$ is perfect, 
\item [$(2)$] $M \subseteq \fkm F$, and 
\item [$(3)$] $\fkm M=\fkm N$, 
\end{enumerate}
then the equality $M^2=NM$ holds. 
\end{Lemma}

\begin{proof}
By assumption $(2)$, $M^2 \subseteq \fkm MF \subseteq NF$. 
Take $x \in M^2$, and write $x=\sum_{i=1}^n h_ic_i, $ $(h_i \in F)$. 
For each $\alpha \in \fkm$, 
$$\alpha x=\sum_{i=1}^n(\alpha h_i)c_i \in \fkm M^2=\fkm N^2 \subseteq N^2$$
by assumption $(3)$. 
By Corollary 2.5, each coefficient $\alpha h_i$ of $c_i$ is in $N$. 
Therefore each $h_i \in N:_F \fkm =M$. 
Hence $x \in NM$. We get the equality $M^2=NM$. 
\end{proof}

The following is the first class of modules of reduction number one, which was 
first discovered by Corso and Polini \cite{CP} when $r=1$. 

\begin{Theorem}
Let $N$ be a submodule of $F$ generated by $c_1, \dots , c_n$, where 
$n=\mu_A(N)$. Assume $n>r$. Put $M=N:_F \fkm$ the socle module of $N$. 
Assume that 
\begin{enumerate}
\item[$(1)$] $\tilde{N}$ is perfect, 
\item[$(2)$] $N \subseteq \fkm F$, and 
\item[$(3)$] $\depth A > 0$. 
\end{enumerate}
Then the equality $M^2=NM$ holds, if $A$ is not a regular local ring. 
\end{Theorem}

\begin{proof}
By assumption $(1)$, the projective dimension of $N$ is finite, because 
the Buchsbaum-Rim complex of the matrix $\tilde{N}$ is acyclic. 
By assumption $(3)$ and $A$ is not regular, we have the equality 
$\fkm N=\fkm M$ (Lemma 3.1). 
Consequently, to prove $M^2=NM$, it is enough to show that 
$M \subseteq \fkm F$. 
Assume $M \nsubseteq \fkm F$ and take $t \in M \backslash \fkm F$. 
Since $t$ is a part of a free basis for $F$, the inclusion $\fkm t 
\hookrightarrow \fkm F$ splits. 
Since $\fkm t \subseteq N \subseteq \fkm F$, the inclusion 
$\fkm t \hookrightarrow N$ also splits. 
Hence the maximal ideal $\fkm$ of $A$ is a direct summand of $N$ 
with finite projective dimension. But this is a contradiction to 
the assumption $A$ is not regular. 
Hence $M \subseteq \fkm F$. Therefore we get the equality  
$M^2=NM$ by Lemma 4.1. 
\end{proof}

To construct the other class of modules of reduction number one, 
we need the following lemma. It is a generalization of 
\cite[Lemma 3.5]{CHV}. 

\begin{Lemma}
Let $(B, \fkn)$ be an Artinian local ring and $W$ a finitely generated 
$B$-module with $\dim_{B/\fkn}(\Soc_B (W))=1$. 
Let $V$ be a submodule of $W$ generated by $x_1, \dots , x_{\ell}$, 
where $\ell=\mu_B(V)$. 
Then, for any non-zero socle $\Delta \in \Soc_B(W)$, there exist 
$b_1, \dots , b_{\ell} \in B$ such that 
$$b_ix_j=\delta_{ij} \Delta$$
for all $1 \leq i, j \leq \ell$ where $\delta_{ij}$ denotes the Kronecker delta. 
\end{Lemma}

\begin{proof}
Fix an integer $1 \leq i \leq \ell$. We look at the following commutative diagram:
$$
\xymatrix{
V \ar[rrr] \ar@{>>}[d] & & & W \ar[r] & E_B(W) \ar[d]^{\varphi_i} \\
V/\fkn V \ar[r]^{p_i} & B/\fkn \ar[r]^{\Delta} & \Delta B \ar[r] & W \ar[r]^{f} & E_B(W), 
}
$$
where $E_B(W)$ is an injective hull of $W$, $\Delta$ is a multiplication map and $p_i(\bar{x_j}) = \delta_{ij}$. 
Since $\dim_{B/\fkn}(\Soc_B (W))=1$, $E_B(W) \cong E_B(B/\fkn)$. Hence $\varphi_i \in \End_B(E_B(B/\fkn)) \cong B$ so that it is a multiplication map for some element $b_i \in B$. 
Since the map $f$ is injective, we have 
$$
\xymatrix{
V \ar[r] \ar@{>>}[d] & W \ar[d]^{b_i} \\
V/\fkn V \ar[r]^{\Delta p_i} & W. 
}
$$
By this commutative diagram, we have $b_i x_j = \delta_{ij} \Delta $
for all $1 \leq j \leq \ell$. 
\end{proof}

Now let me give the other class of modules of reduction number one, which 
was first discovered by Corso, Huneke and Vasconcelos \cite[Theorem 3.7]{CHV}
when $r=1$. 

\begin{Theorem}
Let $N$ be a submodule of $F$ generated by $c_1, \dots , c_n$, where 
$n=\mu_A(N)$. Assume $n>r$. Put $M=N:_F \fkm$ 
the socle module of $N$. 
Assume that 
\begin{enumerate}
\item[$(1)$] $N \subseteq \fkm F$, 
\item[$(2)$] $\ell_A(F/N) < \infty$, and 
\item[$(3)$] $\dim_k(\Soc_A(F/N))=1$. 
\end{enumerate}
Then the equality $M^2=NM$ holds, if $\mu_A(\fkm F/N) \geq 2$. 
\end{Theorem}

\begin{proof}
Let $\ell=\mu_A(\fkm F/N)$ and $\fkm F/N=\sum_{i=1}^{\ell}A \bar{x_i}$, where
$x_i \in \fkm F$. 
Let $I(N)=\Fitt_0(F/N)$ be the $0$-th Fitting ideal of $F/N$ and 
put $B=A/I(N)$. 
Then $I(N)$ is an $\fkm$-primary ideal, because $\ell_A(F/N) < \infty$. 
Let $V=\fkm F/N \subset W=F/N$ and $\bar{\Delta}$ a non-zero socle of $F/N$. 
By Lemma 4.3, there exist $b_1, \dots , b_{\ell} \in B$ such that 
$b_i \bar{x_j}=\delta_{ij} \bar{\Delta}$. 
Let $b_i=\bar{a_i}$, where $a_i \in A$. 
We can write 
$$ \Delta=a_1x_1+n_1=a_2x_2+n_2 $$
for some $n_1, n_2 \in N$, because $\ell \geq 2$. 
Since $M=N+A \Delta$, it is enough to show that 
$\Delta^2 \in NM$. 
Then 
\begin{eqnarray*}
\Delta^2 &=& (a_1x_1+n_1)(a_2x_2+n_2) \\
         &=& (a_1x_2)(a_2x_1)+a_1x_1n_2+a_2x_2n_1+n_1n_2,  
\end{eqnarray*}
so that we have $\Delta^2 \in NM$. Hence the equality $M^2=NM$ holds. 
\end{proof}

\section{The Main Theorem}

We are now ready for the main result in the paper. 

\begin{Theorem}
Let $(A, \fkm)$ be a Cohen-Macaulay local ring of dimension two 
and $N$ a parameter module in a free module $F=A^r$ of rank $r>0$, 
namely, $\ell_A(F/N) < \infty$ and $\mu_A(N)=r+1$. 
Assume that $N \subseteq \fkm F$ and put $M=N:_F \fkm$ the socle 
module of $N$. Then the following three conditions are equivalent. 
\begin{enumerate}
\item[$(1)$] $M^2 \ne NM$. 
\item[$(2)$] $N$ is an integrally closed module in $F$. 
\item[$(3)$] $A$ is a regular local ring, $r=1$ and $N$ is an integrally closed ideal in $A$. 
\end{enumerate}
\end{Theorem}

\begin{proof}
As stated in the introduction, it is known in the case where $r=1$. 

$(2) \Rightarrow (3)$: Since $N$ is a parameter module over a Cohen-Macaulay local ring, $\pd_A N < \infty$. Therefore $A$ is regular by Corollary 3.2. Since $N$ is integrally closed in $F$, $N$ is $\fkm$-full (\cite[Proposition 2.6]{BVa}). Hence we have 
$$r+1=\mu_A(N) \geq \mu_A(\fkm F)=2r, $$ 
by \cite[Corollary 2.7]{BVa}. We get $r=1$. 

$(1) \Rightarrow (3)$: 
Assume $M^2 \ne NM$. 
Since $N$ is a parameter module over a Cohen-Macaulay local ring, 
the matrix $\tilde{N}$ is perfect. So we have 
$A$ is a regular local ring by Theorem 4.2. 
Since $d=2$, we get that $\dim_k(\Soc_A(F/N))=1$. 
Hence we have $\mu_A(\fkm F/N) \leq 1$ by Theorem 4.4. 
Consequently, 
$$\mu_A(\fkm F) = 2r \leq \mu_A(N)+1 =(r+1)+1 =r+2. $$
Therefore $r \leq 2$. 
Assume $r=2$. Then we have the following. 

\smallskip

\noindent
{\bf{Claim}} \begin{enumerate}
\item [$(i)$] $M \subseteq \fkm F $. 
\item [$(ii)$] $\fkm M =\fkm N$. 
\end{enumerate}

\begin{proof}[Proof of Claim]
$(i)$: Assume $M \nsubseteq \fkm F$. Similar to the proof of Theorem 4.2, 
we have that the maximal ideal $\fkm$ is a direct summand of $N$. 
Replacing the generators of $N$ if necessary, 
we may assume that 
$$\tilde{N}=
\left(
\begin{array}{ccc}
x & y & \alpha \\
0 & 0 & \beta  
\end{array}
\right), 
$$
where $\fkm=(x, y)$ and $\alpha, \beta \in \fkm$ because $\mu_A(N)=3$. 
Then $I_2(\tilde{N})=\fkm \beta \subseteq (\beta)$, which is contradict to 
$\ell_A(F/N) < \infty$. Consequently, we have $M \subseteq \fkm F$. 

$(ii)$: Since $r=2$, $$\mu_A(\fkm F)=\mu_A(N)+1=4. $$
So the minimal system of generators of $N$ is a part of a system of 
generators of $\fkm F$. This implies that $N \cap \fkm^2 F=\fkm N$. 
By Claim $(i)$, $\fkm M \subseteq \fkm ^2 F \cap N=\fkm N$. 
Therefore we have $\fkm M=\fkm N$. 
\end{proof}

\noindent
Consequently, we have the equality $M^2=NM$ by Lemma 4.1. This is a contradiction. 
Hence $r=1$. 
\end{proof}

As a direct consequence, we have the following corollaries. 

\begin{Corollary}
Let $(A, \fkm)$ be a two dimensional Cohen-Macaulay local ring and 
$N$ a parameter module in $F=A^r$. Assume that $N \subseteq \fkm F$
and put $M=N:_F \fkm$. Then the equality $M^2=NM$ holds, 
if one of the following assertions is satisfied:
\begin{enumerate}
\item [$(1)$] $A$ is not a regular local ring, 
\item [$(2)$] the rank $r$ is at least two. 
\end{enumerate}
When this is the case, 
the Rees algebra $\Rees(M)$ of $M$ is Cohen-Macaulay. 
\end{Corollary}

\begin{Corollary}
Let $(A, \fkm)$ be a regular local ring of dimension two and $I$ an $\fkm$-primary ideal in $A$. Let $n=\mu_A(I)$. Consider the minimal free 
resolution of $A/I$:
$$
\xymatrix{
F_{\bullet} : 0 \ar[r] & A^{n-1} \ar[r]^{\varphi} & 
A^n \ar[r] & A \ar[r] & A/I \ar[r] & 0. 
}
$$
Taking the $A$-dual of this complex, we have the exact sequence:
$$ 
\xymatrix{
F_{\bullet}^{\ast} : 0 \ar[r] & A \ar[r] & A^n \ar[r]^{\varphi^{\ast}} & 
A^{n-1} \ar[r] & \omega_{A/I} \ar[r] & 0,  
}
$$
where $\omega_{A/I}$ is the canonical module of $A/I$. 
Let $N=\Im \varphi^{\ast}$ and $F=A^{n-1}$. We put $M=N:_F \fkm$ the socle module of $N$. Then the equality $M^2=NM$ holds, if $n \geq 3$. 
When this is the case, the Rees algebra $\Rees(M)$ of $M$ is Cohen-Macaulay. 
\end{Corollary}

When the dimension of $A$ is one, one can easily construct an integrally closed 
parameter module of any rank $r$ as follows:
Let $(A, \fkm)$ be a discrete valuation ring and $N=\fkm \oplus \dots \oplus \fkm$ a direct sum of $r$-copies 
of the maximal ideal $\fkm$. Then $N$ is an integrally closed parameter module in $F=A^r$. Also, one can easily check that $M^2 \ne NM$ where $M$ is the socle module of $N$. So, Theorem 5.1 is not true in the case where the dimension of $A$ is one. However, I don't know whether Theorem 5.1 is true or not in the case where the dimension of $A$ is at least three. At the end of this paper, we see the following simple example of the socle module of a typical parameter module over a $3$-dimensional regular local ring. 

\begin{Example}
Let $(A, \fkm)$ be a regular local ring of dimension three and $F=A^2$ a free module. Let $\{ x, y, z \}$ be a regular system of parameters in $A$. 
Consider a parameter module $N$ in $F$ whose matrix $\tilde{N}$ is given by 
$$
\left(
\begin{array}{cccc}
x & y & z & 0 \\
0 & x & y & z 
\end{array}
\right). 
$$
Put $M=N:_F \fkm$. 
Since the minimal free resolution of $F/N$ is given by the Buchsbaum-Rim complex and $\Hom_A (A/\fkm, F/N) \cong \Tor_3^A(A/\fkm, F/N)$, we have $\dim_{A/\fkm}\Soc_A(F/N)=2$. Then it is easy to see that $M=\fkm \oplus \fkm$ and $M^2=NM$ holds. 
\end{Example}

\section*{Acknowledgments}

The work on this paper is based on the talk at the University of Helsinki in the summer of 2005. 
The author thanks Professor Eero Hyry and the members of his seminar--Doctors Tarmo J\"arvilehto and Lauri Ojala--for their hospitality during the stay. 
He also thanks Professor Shiro Goto for his valuable advice and encouragement during this research. The author was partly supported by the Academy of 
Finland, project 48556. 



\end{document}